\documentclass[11pt]{article}
\usepackage{amsmath,amssymb,amsfonts}
\usepackage[latin1]{inputenc}
\usepackage{epsfig}
\numberwithin{equation}{section}
\newtheorem{thm}{Theorem}[section]

\newtheorem{alem}[thm]{Lemma}
\newtheorem{aprop}[thm]{Proposition}
\newtheorem{acor}[thm]{Corollary}
\newtheorem{arem}[thm]{Remark}
\newenvironment{adem}[1][]%
   {\ \\ {\bf Proof #1. }}%
   {\hfill\mbox{\rule{2 true mm}{3 true mm}}}%\vskip 2 ex\noindent}
   {\ \\ {\bf Example #1. }}%
   {\hfill\mbox{\rule{2 true mm}{3 true mm}}}%\vskip 2 ex\noindent}
\newcommand{\R}{\mathbb{R}}
\renewcommand{\P}{\mathbb{P}}

\newcommand{\E}{{\mathbb E}}

\newcommand{\N}{{\mathbb N}}

 \title{On the long time behaviour of stochastic vortices systems}

 \author{J.Fontbona\thanks{DIM-CMM, UMI(2807) UCHILE-CNRS, Universidad de Chile, Casilla 170-3, Correo 3,
Santiago-Chile, e-mail:fontbona@dim.uchile.cl. Partially supported by
Fondecyt Grant 1110923, 
BASAL-Conicyt Center for Mathematical Modeling and Millennium Nucleus  NC130062}\and B.Jourdain\thanks{Université Paris-Est, CERMICS, 6-8 av Blaise Pascal, Cit\'e
     Descartes, Champs sur Marne, 77455 Marne-la-Vall\'ee Cedex 2, France, and INRIA Paris-Rocquencourt - MATHRISK . 
     e-mail : jourdain@cermics.enpc.fr. Supported by the French National
Research Agency (ANR) under the program ANR-12-BLAN Stab}
    }

\begin{document}
 \maketitle
 \begin{abstract}
In this paper, we are interested in the long-time behaviour of stochastic systems of $n$ interacting vortices: the position in $\R^2$ of each vortex evolves according to a Brownian motion and a drift summing the influences of the other vortices computed through the Biot and Savart kernel and multiplied by their respective vorticities. For fixed $n$, we perform the rescalings of time and space used successfully by Gallay and Wayne \cite{gawa} to study the long-time behaviour of the vorticity formulation of the two dimensional incompressible Navier-Stokes equation, which is the limit as $n\to\infty$ of the weighted empirical measure of the system under mean-field interaction. When all the vorticities share the same sign, the $2n$-dimensional process of the rescaled positions of the vortices is shown to converge exponentially fast as time goes to infinity to some invariant measure which turns out to be Gaussian if all the vorticities are equal. In the particular case $n=2$ of two vortices, we prove exponential convergence in law of the $4$-dimensional process to an explicit random variable, whatever the choice of the two vorticities. We show that this limit law is not Gaussian when the two vorticities are not equal.


\end{abstract}
\section{Introduction}
In this work, we are interested in stochastic systems of
interacting vortices :
\begin{equation}
 X^{i}_t=X^i_0+\sqrt{2\nu}\;W^i_t+\int_0^t\sum_{\stackrel{j=1}{j\neq
    i}}^n a_j K(X^{i}_s-X^{j}_s)ds,\;1\leq i\leq n
\label{systpart}
\end{equation}
 where $K:x=(x_1,x_2)\in\R^2\rightarrow \frac{x^\perp}{2\pi|x|^2}$ ($x^\perp=(-x_2,x_1)$) denotes the Biot and Savart kernel,
$(X^i_0)_{i\geq 1}$ are two-dimensional random vectors independent from the sequence $(W^i)_{i\geq
1}$ of independent two-dimensional Brownian motions. For $i\in\{1,\hdots,n\}$, the real number $a_i$ is the intensity or vorticity of the $i$-th vortex. The Biot and
Savart kernel $K$ is singular at the origin but locally bounded and
Lipschitz continuous on $\R^2\setminus\{0\}$. Under the assumption $\P(\exists i\neq j\mbox{ s.t. }X^i_0=X^j_0)=0$ that will be made throughout the paper, existence and uniqueness
results for this $2n$-dimensional stochastic differential equation  are
given in \cite{Tak,O1,O2,FM}. Moreover, it is shown in \cite{O1}  that for $t>0$, the random vector $X_t=(X_t^{1},\dots,X^{n}_t)$ has a density
$\rho_t(x)$.

The relation between System \eqref{systpart}  and the vorticity formulation
\begin{equation}
   \partial_tw(t,x)=\nu \Delta w(t,x)-\nabla.\left(w(t,x)\int_{\R^2}K(x-y)w(t,y)dy\right)\label{nsvort}
\end{equation} of the two dimensional incompressible Navier-Stokes equation is well known and has been studied by several authors. It arises through the propagation of chaos property, stating  in particular that after a suitable normalization of the vortex intensities $a_i$, some weighted empirical measure of the system \eqref{systpart}  or  some variant of it  involving a regularized version of the kernel $K$, converges in law as $n$ tends to $\infty$ to $w_t$, see  \cite{MaPu,Mel1}  for the regularized case and \cite{Oprop,FHM}  that deal with the true Biot and Savart kernel.

On the other hand, the long time behavior for the two dimensional vortex equation has been successfully studied in Gallay and Wayne \cite{gawa}, who established the strong convergence of a time-space rescaled  version of its solution $w_t$ as $t\to \infty$ to a gaussian density with total mass given by the initial circulation $\int_{\R^2} w_0(x)dx$. 

Motivated by that result, our goal  in this paper  is to explore some of the asymptotic properties of System  \eqref{systpart} as time tends to $\infty$. To that end, following the scaling introduced in \cite{gawa},
we define $Z^i_t=e^{-\frac{t}{2}}X^{i}_{e^t-1}$. The process $(Z_t=(Z^1_t,\hdots,Z^n_t))_{t\geq 0}$ solves \begin{equation}
Z^{i}_t 
=X^i_0+\sqrt{2\nu}B^i_t+\int_0^t\bigg(\sum_{\stackrel{j=1}{j\neq
    i}}^na_jK(Z^{i}_s-Z^{j}_s)-\frac{Z^{i}_s}{2}\bigg)ds,\;1\leq
    i\leq n\label{resys}\end{equation} where
$\left(B^i_t=\int_0^{e^t-1}\frac{dW^i_s}{\sqrt{s+1}}\right)_{i\geq
1}$ are independent two dimensional Brownian motions according to
the Dambin-Dubins-Schwarz theorem.
This consequence of the homogeneity of the Biot and Savart kernel can easily be checked by computing $\frac{X^{i}_u}{\sqrt{u+1}}$ by Itô's formula and then choosing $u=e^t-1$.
For $t>0$, $Z_t$ admits the density \begin{equation}
   p_t(z)=e^{nt}\rho_{e^t-1}(ze^{t/2})\label{lienpq}
\end{equation} with respect to the Lebesgue measure on
$\R^{2n}$. Moreover, the Fokker-Planck equation
\begin{equation}\label{FPZ}
\partial_t p_t(z)=\nu \triangle p_t(z) - \sum_{i=1}^n \bigg(\sum_{\stackrel{j=1}{j\neq
    i}}^na_jK(z^i-z^j)\bigg)\cdot \nabla_{z^i} p_t(z) + \frac{1}{2}\nabla_{z}. (z p_t(z)), \ t>0, z\in
    \mathbb{R}^{2n}.
\end{equation}
holds in the weak sense.

The outline of the paper and the results that we obtain are the following. In Section \ref{constsign}, we  consider System \eqref{resys}  with vorticities $a_i$ of constant sign, and show in Subsection \ref{Lyapconstsign}  by using  Lyapunov-Foster-Meyn-Tweedie techniques the exponential convergence in total variation norm to some invariant law with a positive density with respect to the Lebesgue measure. In Subsection \ref{equalvort}, we study the more restrictive case of particles having equal vorticities. In that case, we show that the invariant measure is the standard $2n$ dimensional normal law scaled by the viscosity coefficient. Using the Logarithmic Sobolev inequality  satisfied by this invariant measure, we deduce  that the relative entropy of the law of the particle positions with respect to that Gaussian law goes exponentially fast to $0$, at an explicit rate independent of $n$. Although the techniques employed in the proof of these two results are standard, the singularity of the drift requires some specific adaptations by regularization. These   results are also translated into the original time-space scale.   We have not been able to study by similar techniques the long-time behaviour when the vorticities are allowed to have different signs. In order to gain some insight on the difficulties that the long time behavior raises in that case,  we consider in Section \ref{toy} the particular case of $n=2$ vortices with arbitrary intensities. 
We completely describe the equilibrium law in $\R^4$  in  terms of the stationary solution of  some related stochastic differential equation. In particular, we prove that, although some linear combinations of the two particles positions are Gaussian under the stationary measure,  this stationary measure is not jointly Gaussian unless the intensities are equal. We moreover show exponential convergence to this limiting law in suitable Wasserstein distances. The argument for long-time convergence, which is based on coupling and time reversal techniques, is original to our knowledge and could be of interest in different contexts.

\subsection*{Notation}
\begin{itemize}
   \item For $d\in{\mathbb N}^*$ and $\alpha\geq 1$, let $W_\alpha$ denote the
Wasserstein metric on the space of probability measures on $\R^d$
defined by
$$W_\alpha(\mu,\nu)=\inf_{\rho\tiny{<\begin{array}{l}\mu
   \\\nu
\end{array}}}\left(\int_{\R^{d}\times\R^d}|x-y|^\alpha\rho(dx,dy)\right)^{1/\alpha},$$
where the infimum is computed over all measures $\rho$ on
$\R^d\times\R^d$ with first marginal equal to $\mu$ and second marginal
equal to $\nu$ and $|x-y|$ denotes the Euclidean norm of $x-y\in\R^d$.
\item To deal with the singularity of the Biot and Savart kernel, we construct smooth approximations of this kernel which
 coincide with it away from $0$. Let $\varphi$ be a smooth function on $\R_+$ such that
\begin{equation*}
  \varphi(r)=\begin{cases}\frac{1}{2}\mbox{ for $r\leq\frac{1}{2}$}
    \\r\mbox{ for $r\geq 1$}
\end{cases},
\end{equation*}
and for $\varepsilon>0$, $\varphi_\varepsilon(r)=\varepsilon\varphi(r/\varepsilon)$. The
kernel $K_\varepsilon(z)\stackrel{\rm def}{=}\frac{1}{2\pi}\nabla^\perp_z
\ln(\varphi_\varepsilon(|z|))$ coincides with $K$ for $|z|\geq\varepsilon$ and is
divergence-free, globally Lipschitz continuous and bounded.\end{itemize}

\section{Case of vorticities $(a_i)_{1\leq i\leq n}$ with constant sign}\label{constsign}
In the present section, we assume that either $\forall i\in\{1,\hdots,n\}$, $a_i>0$ or $\forall i\in\{1,\hdots,n\}$, $a_i<0$.
Let ${\cal Z}=\{z=(z^1,\hdots,z^n)\in\R^{2n}:\;\forall i\neq j,\;z^i\neq z^j\}$.
For any $z\in{\cal Z}$, the unique strong solution $(Z^z_t)$ of the SDE \eqref{resys} starting from $(Z^1_0,\hdots,Z^n_0)=z$ is such that $\tau^z=\inf\{t\geq 0:\exists i\neq j,\;Z^{i}_t=Z^{j}_t\}$ is a.s. infinite.
Since the process $(\sqrt{1+t}(Z^{1}_{\ln(1+t)},\hdots,Z^{n}_{\ln(1+t)}))_{t\geq 0}$ solves 
\eqref{systpart}, this is in particular a consequence of \cite{Tak}.

\subsection{Convergence to equilibrium}\label{Lyapconstsign} 
\begin{aprop}\label{experg}
The SDE \eqref{resys} admits a unique invariant probability measure. This measure admits a positive density $p_\infty$ with respect to the Lebesgue measure and there exist constants $C,\lambda\in (0,+\infty)$ depending on $\nu$, $n$ and $(a_1,\cdots,a_n)$ such that
$$\forall z\in{\cal Z},\;\sup_{f:\forall x,|f(x)|\leq 1+V(x)}\left|\E[f(Z^z_t)]-\int_{\R^{2n}}f(x)p_\infty(x)dx\right|\leq C e^{-\lambda t}(1+V(z)),$$
where 
\begin{equation}
   V(z)=\sum_{i=1}^n |a_i||z_i|^2.\label{lyapufunct}
\end{equation} 
\end{aprop}According to \eqref{lienpq}, the density
of the original particle system $(X^1_s,\hdots,X^n_s)$ is
\begin{equation}\rho_s(x)=(1+s)^{-n}p_{\ln(1+s)}\left(\frac{x}{\sqrt{1+s}}\right) \label{scaldens}.\end{equation}
We deduce that
\begin{acor}
Let $x=(x^1,\hdots,x^n)\in{\cal Z}$. Denoting by $X^x_t$ the solution of \eqref{systpart} starting from $(X^1_0,\hdots,X^n_0)=(x^1,\hdots,x^n)$, one has 
$$\sup_{f:\forall y,|f(y)|\leq 1+V(y/\sqrt{1+t})}\left|\E[f(X^x_t)-\int_{\R^{2n}}\frac{f(y)}{(1+t)^n}p_\infty\left(\frac{y}{\sqrt{1+t}}\right)dy\right|\leq \frac{C}{(1+t)^\lambda}\left(1+V\left(\frac{x}{\sqrt{1+t}}\right)\right).$$ 
\end{acor}
All the statements but the existence of the positive density $p_\infty$ are a consequence of Theorem 6.1 in \cite{MT} which supposes a Lyapunov condition and that all compact sets  of ${\cal Z}$ are petite sets for some skeleton chain. Let us first check that the function $V$ defined by \eqref{lyapufunct} is a Lyapunov function. Denoting by $$L=\nu\Delta_z-\frac{z}{2}.\nabla_z+\sum_{i=1}^n\sum_{\stackrel{j=1}{j\neq i}}^n a_j K(z^i-z^j).\nabla_{z^i}$$ the infinitesimal generator associated with \eqref{resys}, one has
\begin{align*}
   LV(z)&=4\nu\sum_{i=1}^n|a_i|-V(z)+2\sum_{j\neq i} a_ia_j K(z^i-z^j).z^i=4\nu\sum_{i=1}^n|a_i|-V(z)+\sum_{j\neq i} a_ia_j \frac{(z^i-z^j)^\perp.(z^j-z^i)}{2\pi|z^i-z^j|^2}\\&=4\nu\sum_{i=1}^n|a_i|-V(z)\end{align*}
where we used the oddness of the Biot and Savart kernel $K$ for the second equality. Hence $$LV(z)\leq \left(4\nu\sum_{i=1}^n|a_i|\right)1_{\{\sum_{i=1}^n |a_i||z_i|^2\leq 8\nu \sum_{i=1}^n|a_i|\}}-\frac{V(z)}{2}.$$
and condition (CD3) in Theorem 6.1 \cite{MT} is satisfied. So it only remains to check that all compact sets of ${\cal Z}$ are petite sets for some skeleton chain to apply this theorem. This is a consequence of Proposition 6.2.8 \cite{MTbook} and of the next Lemma which implies that any skeleton chain is Lebesgue measure-irreducible and that the invariant probability measure admits a positive density.
\begin{alem}\label{feller}
   The semi-group $(P_t)_{t\geq 0}$ defined by $P_tf(z)=\E[(f(Z^z_t))]$ is Feller and for any $z\in{\cal Z}$ and any $t>0$, $Z^z_t$ admits a positive density with respect to the Lebesgue measure.
\end{alem}
\begin{adem}
Let us check by probabilistic estimations that the semi-group $(P_t)_{t\geq 0}$ is Feller and first that for $f:\R^{2n}\to\R$ bounded and going to $0$ at infinity, so does $P_tf$.\\
Let $R^z_t=\sum_{i=1}^n|a_i||Z^{z,i}_t|^2$. By Itô's formula,
\begin{align}
   dR^z_t&=2\sqrt{2\nu}\sum_{i=1}^n |a_i|Z^{z,i}_t.dB^i_t-R^z_tdt+4\nu\sum_{i=1}^n |a_i|dt.\label{dynr}\\
d\frac{1}{1+R^z_t}&=-\frac{2\sqrt{2\nu}}{(1+R^z_t)^2}\sum_{i=1}^n |a_i|Z^{z,i}_t.dB^i_t+\frac{R^z_t}{(1+R^z_t)^2}dt-\frac{4\nu\sum_{i=1}^n |a_i|}{(1+R^z_t)^2}dt+\frac{8\nu}{(1+R^z_t)^3}\sum_{i=1}^n a_i^2|Z^{z,i}_t|^2dt.\notag\end{align}
The stochastic integral $\left(\int_0^t\frac{1}{(1+R^z_s)^2}\sum_{i=1}^n |a_i|Z^{z,i}_s.dB^i_s\right)_{t\geq 0}$ is a locally in time square integrable martingale. Denoting by $\bar{a}=\max_{1\leq i\leq n}|a_i|$, remarking that $\sum_{i=1}^n a_i^2|Z^{z,i}_t|^2\leq \bar{a}R^z_t$ and taking expectations, we deduce that 
$$\frac{d}{dt}\E\left[\frac{1}{1+R^z_t}\right]\leq \E\left[\frac{R^z_t}{(1+R^z_t)^2}+\frac{8\nu\bar{a}R^z_t}{(1+R^z_t)^3}\right]\leq (1+8\nu\bar{a})\E\left[\frac{1}{1+R^z_t}\right].$$
Therefore
\begin{equation}
\E\left[\frac{1}{1+\bar{a}|Z^z_t|^2}\right]\leq \E\left[\frac{1}{1+R^z_t}\right]\leq \frac{1}{1+\sum_{i=1}^n|a_i||z^i|^2}e^{ (1+8\nu\bar{a})t}\leq \frac{1}{1+|z|^2\min_{1\leq i\leq n}|a_i|}e^{ (1+8\nu\bar{a})t}\label{momneg}. 
\end{equation}
For $r\in(0,+\infty)$, since by Markov inequality, $\P(|Z^z_t|\leq r)\leq\E\left[\frac{1+\bar{a}r^2}{1+\bar{a}|Z^z_t|^2}\right]$, we deduce that $\lim_{|z|\to+\infty}\P(|Z^z_t|\leq r)=0$. As a consequence, for any bounded function $f$ on $\R^{2n}$ going to $0$ at infinity, $\lim_{|z|\to\infty}P_tf(z)=0$.

Let $f:\R^{2n}\to\R$ be continuous and bounded. To check that $z\in{\cal Z}\mapsto P_tf(z)$ is continuous, we introduce for $\varepsilon>0$  a bounded and smooth kernel $K_\varepsilon:\R^2\to\R^2$ which coincides with the Biot and Savart kernel $K$ on $\{x\in\R^2:|x|\geq \varepsilon\}$ and define $Z^{z,\varepsilon}=(Z^{z,\varepsilon,1},\hdots,Z^{z,\varepsilon,1})$ as the solution of the SDE 
\begin{equation}
Z^{z,\varepsilon,i}_t 
=z^i+\sqrt{2\nu}B^i_t+\int_0^t\bigg(\sum_{\stackrel{j=1}{j\neq
    i}}^na_jK_\varepsilon(Z^{z,\varepsilon,i}_s-Z^{z,\varepsilon,j}_s)-\frac{Z^{z,\varepsilon,i}_s}{2}\bigg)ds,\;1\leq
    i\leq n\label{regsys}.\end{equation}Since $K_\varepsilon$ is Lipschitz continuous, one easily checks that $|Z^{z,\varepsilon}_t-Z^{z',\varepsilon}_t|\leq |z-z'|e^{C_\varepsilon t}$ for some deterministic finite constant $C_\varepsilon$ depending on $\varepsilon$ but not on $z,z'\in\R^{2n}$. Hence $z\mapsto \E[f(Z^{z,\varepsilon}_t)]$ is continuous. To deduce the continuity of $z\mapsto \E[f(Z^{z}_t)]$, we are going to control the stopping times 
$$\tau^{z,\varepsilon}=\inf\{t\geq 0:\exists i\neq j:|Z^{z,i}_t-Z^{z,j}_t|\leq\varepsilon\}.$$
Notice that for $z\in{\cal Z}$, the processes $Z^z$ and $Z^{z,\varepsilon}$ coincide on $[0,\tau^{z,\varepsilon}]$.
Therefore,
\begin{equation}
\forall z,z'\in{\cal Z},\;|P_tf(z)-P_tf(z')|\leq |\E[f(Z^{z,\varepsilon}_t)-f(Z^{z',\varepsilon}_t)]|+2\|f\|_\infty\left(\P(\tau^{z,\epsilon}<t)+\P(\tau^{z',\epsilon}<t)\right).\label{contsg}
\end{equation}
By continuity of the paths of $Z^z$, for $z\in{\cal Z}$, $\lim_{\varepsilon\to 0}\tau^{z,\varepsilon}=\tau^z=+\infty$ a.s. which implies that $\lim_{\varepsilon\to 0}\P(\tau^{z,\epsilon}<t)=0$. But we need some uniformity in the starting point $z$ to deduce that $z\in{\cal Z}\mapsto P_tf(z)$ is continuous. The following computations  inspired from \cite{Tak} improve the result therein into a quantitative estimate. By Itô's formula and since
\begin{align*}\sum_{i\neq j}&a_ia_j\frac{Z^{z,i}_t-Z^{z,j}_t}{|Z^{z,i}_t-Z^{z,j}_t|^2}.\left(\sum_{\stackrel{k=1}{k\neq i}}^na_k\frac{(Z^{z,i}_t-Z^{z,k}_t)^\perp}{|Z^{z,i}_t-Z^{z,k}_t|^2}-\sum_{\stackrel{l=1}{l\neq j}}^na_l\frac{(Z^{z,j}_t-Z^{z,l}_t)^\perp}{|Z^{z,j}_t-Z^{z,l}_t|^2}\right)\\
&=2\sum_{i=1}^na_i\left(\sum_{\stackrel{j=1}{j\neq i}}^na_j\frac{Z^{z,i}_t-Z^{z,j}_t}{|Z^{z,i}_t-Z^{z,j}_t|^2}\right).\left(\sum_{\stackrel{k=1}{k\neq i}}^na_k\frac{Z^{z,i}_t-Z^{z,k}_t}{|Z^{z,i}_t-Z^{z,k}_t|^2}\right)^\perp=0,
\end{align*}
one obtains
\begin{align*}
   d\sum_{i\neq j}a_ia_j\ln|Z^{z,i}_t-Z^{z,j}_t|=&\sqrt{2\nu}\sum_{i\neq j}a_ia_j\frac{Z^{z,i}_t-Z^{z,j}_t}{|Z^{z,i}_t-Z^{z,j}_t|^2}.(dB^i_t-dB^j_t)-\frac{1}{2}\sum_{i\neq j}a_ia_jdt
\end{align*}
Hence $$\E\left[\sum_{i\neq j}a_ia_j\ln|Z^{z,i}_{\tau^{z,\varepsilon}\wedge t}-Z^{z,j}_{\tau^{z,\varepsilon}\wedge t}|\right]=\sum_{i\neq j}a_ia_j\ln|z^i-z^j|-\frac{1}{2}\sum_{i\neq j}a_ia_j\E[\tau^{z,\varepsilon}\wedge t]\geq \sum_{i\neq j}a_ia_j\left(\ln|z^i-z^j|-\frac{t}{2}\right).$$
Since $\forall x>0,\;\ln(x)<x$, for $\varepsilon\in(0,1)$, the left-hand-side is not greater than
\begin{align*}
   \min_{i\neq j}a_ia_j\ln(\varepsilon)\P(\tau^{z,\varepsilon}\leq t)&+\sum_{i\neq j}a_ia_j\E\left[\sup_{s\in[0,t]}|Z^{z,i}_{s}-Z^{z,j}_s|\right]\\
&\leq \min_{i\neq j}a_ia_j\ln(\varepsilon)\P(\tau^{z,\varepsilon}\leq t)+2\sum_{i=1}^na_i\left(\sum_{\stackrel{j=1}{j\neq i}}^na_j\right)\E\left[\sup_{s\in[0,t]}|Z^{z,i}_{s}|\right]\\&\leq \min_{i\neq j}a_ia_j\ln(\varepsilon)\P(\tau^{z,\varepsilon}\leq t)+2\sum_{i=1}^n\sqrt{|a_i|}\left(\sum_{\stackrel{j=1}{j\neq i}}^n|a_j|\right)\E^{1/4}\left[\sup_{s\in[0,t]}(R^z_s)^2\right].
\end{align*}
Hence
\begin{equation}
  \P(\tau^{z,\varepsilon}\leq t)\leq  \frac{\sum_{i\neq j}a_ia_j\left(\frac{t}{2}-\ln|z^i-z^j|\right)+2\sum_{i=1}^n\sqrt{|a_i|}\left(\sum_{j\neq i}|a_j|\right)\E^{1/4}\left[\sup_{s\in[0,t]}(R^z_s)^2\right]}{\min_{i\neq j}a_ia_j\ln(1/\varepsilon)}\label{majoprobarenc}
\end{equation}
By a standard localisation procedure, one checks that the stochastic integral the differential of which appears in the right-hand-side of \eqref{dynr} is a martingale and that
\begin{equation*}
\E[R^z_t]=e^{-t}\sum_{i=1}^n|a_i||z^i|^2+(1-e^{-t})4\nu\sum_{i=1}^n |a_i|.
\end{equation*}
By \eqref{dynr}, $R^z_s\leq R^z_0+4\nu\sum_{i=1}^n |a_i|t+2\sqrt{2\nu}\int_0^s\sum_{i=1}^n|a_i|Z^{z,n,i}_u.dB^i_u$ and by Doob's inequality, 
\begin{align*}
   \E&\left[\sup_{s\in[0,t]}(R^z_s)^2\right]\leq 2\left((\sum_{i=1}^n|a_i|(|z^i|^2+4\nu t))^2+32\nu\int_0^t\sum_{i=1}^na_i^2\E[|Z^{z,i}_s|^2]ds\right)\\
&\leq 2\left((\sum_{i=1}^n|a_i|(|z^i|^2+4\nu t))^2+32\nu\bar{a}\left(\sum_{i=1}^n|a_i||z^i|^2(1-e^{-t})+4\nu\sum_{i=1}^n|a_i|(t+e^{-t}-1)\right)\right).
\end{align*}
Plugging this estimation in \eqref{majoprobarenc}, we deduce that for $z\in{\cal Z}$ and $\alpha>0$ small enough so that the ball $B(z,\alpha)$ centered in $z$ with radius $\alpha$ is contained in ${\cal Z}$, one has $\lim_{\varepsilon \to 0}\sup_{z'\in B(z,\alpha)}\P(\tau^{z',\varepsilon}\leq t)=0$. With \eqref{contsg} and the continuity of $z\mapsto \E[f(Z^{z,\varepsilon}_t)]$ for fixed $\varepsilon>0$, we conclude that $z\mapsto P_tf(z)$ is continuous.\par

Last, by Theorem 2 and Example 2 \cite{O1}, for $t>0$, $\sqrt{1+t}Z^{z}_{\ln(1+t)}$ admits a density satisfying some Gaussian lower bound. This implies that for $t>0$, $Z^z_t$ admits a positive density. Notice that one could also deduce the Feller property of the semi-group from the estimations of the fundamental solution of $\frac{\partial}{\partial t}-L$ obtained by partial differential equation's techniques in that paper but we preferred to give a probabilistic argument.

\end{adem}

\subsection{Case of particles with equal vorticity}\label{equalvort}

In the present subsection, we assume the existence of $a\in\R^*$ such that $a_i=a$ for all
$i\in\{1,\hdots,n\}$.
For instance, when $w_0$ is a non-negative (resp. non-positive) initial vorticity density
on $\R^2$ with positive and finite total mass $\|w_0\|_1$, it is natural
to choose the initial positions 
$X^i_0$ i.i.d. according to the probability density $\frac{|w_0|}{\|w_0\|_1}$ and
$a=\frac{\|w_0\|_1}{n}$ (resp. $a=-\frac{\|w_0\|_1}{n}$).
In this  situation,  the invariant density turns out to be Gaussian and we explicit the exponential rate of convergence to equilibrium. In fact, both the invariant density and the rate of convergence are the same as for the Ornstein Uhlenbeck dynamics given by \eqref{resys} in the case without interaction: $a_i=0$ for all
$i\in\{1,\hdots,n\}$.

\begin{aprop}Assume the existence of $a\in\R^*$ such that $a_i=a$ for all
$i\in\{1,\hdots,n\}$.
 The density
   $p_\infty(z^1,\hdots,z^n)=\frac{1}{(4\pi\nu)^{n}}e^{-\frac{1}{4\nu}\sum_{i=1}^n|z^i|^2}$of the normal law ${\cal N}_{2n}(0,2\nu I_{2n})$ is invariant for the SDE
    \eqref{resys}. Moreover, if $(X^{1}_0,\hdots,X^n_0)$ has a density $p_0$ with respect
to the Lebesgue measure on $\R^{2n}$ such that the relative entropy $\int_{\R^{2n}}p_0\ln\left(\frac{p_0}{p_\infty}\right)$ is finite. Then, one has
$$\forall t\geq 0,\;\int_{\R^{2n}}
   p_t\ln\left(\frac{p_t}{p_\infty}\right)\leq e^{-t}\int_{\R^{2n}}
   p_0\ln\left(\frac{p_0}{p_\infty}\right).$$
\label{convdensnorm}\end{aprop}
The above rate $1$ of exponential convergence is the same as for the Ornstein-Uhlenbeck process obtained when $a=0$ and depends neither on $n$ nor on $\nu$.\\
 Since, by \eqref{scaldens}, the density $\rho_s$ of the orginal particle system $(X^1_s,\hdots,X^n_s)$ is such that $$\int_{\R^{2n}}
   \rho_s\ln\left(\frac{\rho_s}{(1+s)^{-n}p_\infty\left(\frac{.}{\sqrt{1+s}}\right)}\right)=\int_{\R^{2n}}p_{\ln(1+s)}\ln\left(\frac{p_{\ln(1+s)}}{p_\infty}\right),$$ one easily deduces its asymptotic behaviour as $s\to +\infty$.
\begin{acor}\label{corconvdensnorm}Assume that $(X^{1}_0,\hdots,X^n_0)$ has a density $p_0$ with respect
to the Lebesgue measure on $\R^{2n}$ such that the relative entropy $\int_{\R^{2n}}p_0\ln\left(\frac{p_0}{p_\infty}\right)$ is finite. Then
   $$\forall s\geq 0,\;\int_{\R^{2n}}
   \rho_s\ln\left(\frac{\rho_s}{(1+s)^{-n}p_\infty\left(\frac{.}{\sqrt{1+s}}\right)}\right)\leq \frac{1}{1+s}\int_{\R^{2n}}
   \rho_0\ln\left(\frac{\rho_0}{p_\infty}\right)=\frac{1}{1+s}\int_{\R^{2n}}
   p_0\ln\left(\frac{p_0}{p_\infty}\right).$$
\end{acor}
  
\begin{arem}Assume that $(X^{1}_0,\hdots,X^n_0)$  is such that $\E(\sum_{i=1}^n|X^i_0|^2)<+\infty$ and let $s>0$. We have $\E
\left(\sum_{i=1}^n|X^i_s|^2\right)=\E(\sum_{i=1}^n|X^i_0|^2)+2n\nu s$. Moreover, by Theorem 2 and Example 2 \cite{O1}, $(X^1_s,\hdots,X^n_s)$ admits a density $\rho_s$ bounded by $C_ns^{-n}$ so that
$$\int_{\R^{2n}}
   \rho_s\ln\left(\frac{\rho_s}{(1+s)^{-n}p_\infty\left(\frac{.}{\sqrt{1+s}}\right)}\right)\leq \ln(C_n)+n\ln\left(1+\frac{1}{s}\right)+n\ln(4\pi\nu)+\frac{1}{4\pi\nu(1+s)}\E
\left(\sum_{i=1}^n|X^i_s|^2\right)<+\infty.$$
By Proposition \ref{convdensnorm} and Corollary \ref{corconvdensnorm},
\begin{align*}
   \forall t\geq \ln(1+s),\;&\int_{\R^{2n}}
   p_t\ln\left(\frac{p_t}{p_\infty}\right)\leq (1+s)e^{-t}\int_{\R^{2n}}
   \rho_s\ln\left(\frac{\rho_s}{(1+s)^{-n}p_\infty\left(\frac{.}{\sqrt{1+s}}\right)}\right)\\
\forall t\geq s,\;&\int_{\R^{2n}}
   \rho_t\ln\left(\frac{\rho_t}{(1+t)^{-n}p_\infty\left(\frac{.}{\sqrt{1+t}}\right)}\right)\leq \frac{1+s}{1+t}\int_{\R^{2n}}
   \rho_s\ln\left(\frac{\rho_s}{(1+s)^{-n}p_\infty\left(\frac{.}{\sqrt{1+s}}\right)}\right).
\end{align*}
Notice that $\inf_{s>0}(1+s)\left(\ln(C_n)+n\ln\left(1+\frac{1}{s}\right)+n\ln(4\pi\nu)\right)+\frac{\E
\left(\sum_{i=1}^n|X^i_0|^2\right)+2n\nu s}{4\pi\nu}$ is attained at a unique point $s_\star>0$.
\end{arem}
\begin{arem}
Asume that the initial conditions $(X^i_0)_{1\leq i\leq n}$ are i.i.d. according to some density $p^1_0$ on $\R^2$ such that $\int_{\R^2}p^1_0\ln\left(\frac{p^1_0}{p^1_\infty}\right)<+\infty$ where $p^1_\infty(z^1)=\frac{1}{4\pi\nu}e^{-\frac{1}{4\nu}|z^1|^2}$. Then, for all $t\geq 0$, $(Z^1_t,\hdots,Z^n_t)$ is exchangeable and, denoting by $p^{1,n}_t$ the common density of the random vectors $Z^i_t$, we have
\begin{align*}
   ne^{-t}\int_{\R^2}p^1_0\ln\left(\frac{p^1_0}{p^1_\infty}\right)&=e^{-t}\int_{\R^{2n}}p_0\ln\left(\frac{p_0}{p_\infty}\right)\geq \int_{\R^{2n}}p_t\ln\left(\frac{p_t}{p_\infty}\right)dz\\&= \int_{\R^{2n}}p_t(z)\ln\left(\frac{p_t(z)}{\prod_{i=1}^np^{1,n}_t(z^i)}\right)dz+\int_{\R^{2n}}p_t(z)\sum_{i=1}^n\ln\left(\frac{p^{1,n}_t(z^i)}{p^1_\infty(z^i)}\right)dz\\
&\geq 0+n\int_{\R^2}p^1_t\ln\left(\frac{p^{1,n}_t}{p^1_\infty}\right).
\end{align*}
We have $p^{1,n}_t(.)=e^t\rho^{1,n}_{e^t-1}(e^{\frac{t}{2}}.)$, where $\rho^{1,n}_{s}$ denotes the common density of the vectors $X^{i}_s$. By the transport inequality satisfied by the Gaussian density $p^1_\infty$, $\int_{\R^2}p^1_0\ln\left(\frac{p^1_0}{p^1_\infty}\right)<+\infty$ implies $\int_{\R^2}p^1_0(x)(|\ln(p^1_0(x))|+|x|^2)dx<+\infty$. Hence, by \cite{FHM},  $\rho^{1,n}_{s}(x)dx$ converges weakly to $w(s,x)dx$ as $n\to\infty$ with $w(s,x)$ denoting the solution of the vorticity formulation \eqref{nsvort} of the two dimensional incompressible Navier-Stokes equation starting from $w(0,x)=p^1_0(x)$. By the lower semi-continuity of the relative entropy with respect to the weak convergence, we conclude that
$$\forall s\geq 0,\;\int_{\R^2}w(s,.)\ln\left(\frac{w(s,.)}{(1+s)^{-1}p^1_\infty\left(\frac{.}{\sqrt{1+s}}\right)}\right)\leq \frac{1}{1+s}\int_{\R^2}p^1_0\ln\left(\frac{p^{1}_0}{p^1_\infty}\right),$$
which can also be deduced from the proof of Lemma 3.2 \cite{gawa}.
\end{arem}
We shall first give formal arguments for Proposition \ref{convdensnorm}, which will be made rigorous by replacing $K$ with the regularized kernel $K_{\varepsilon}$. The results for the system \eqref{resys} will then be justified by suitable passages to the limit as $\varepsilon\to 0$.

Let us check that $p_\infty$ solves the stationary version of the
Fokker-Planck equation \eqref{FPZ}. Since for
$i\in\{1,\hdots,n\}$, $\nu\nabla_{z^i}
   p_\infty+\frac{z^i}{2}p_\infty=0$, the result follows from
\begin{align*}
\sum_{i=1}^n \bigg(\sum_{\stackrel{j=1}{j\neq
    i}}^nK(z^i-z^j)\bigg)\cdot \nabla_{z^i} p_\infty &=
   -\frac{p_\infty}{2\nu} \sum_{i\neq j}K(z^i-z^j).z^i=\frac{p_\infty}{2\nu} \sum_{i\neq j}K(z^i-z^j).\frac{z^j-z^i}{2}\\
   &=\frac{p_\infty}{8\pi \nu} \sum_{i\neq j}\frac{(z^i-z^j)^\perp.(z^j-z^i)}{|z^i-z^j|^2}=0,
\end{align*}
where we used the oddness of the Biot and Savart kernel $K$ for the
second equality.
\begin{arem}\label{noninv}
 When the vorticities $a_i$ of the
     particles are different, $p_\infty$ is no longer a solution of the stationary Fokker-Planck equation. Indeed, $L^*p_\infty(z)=\frac{p_\infty(z)}{2\nu}\sum_{i\neq j}a_jK(z^i-z^j).z^i=\frac{p_\infty(z)}{2\nu}\sum_{i\neq j}(a_j-a_i)\frac{(z^i-z^j)^\perp.(z^i+z^j)}{8\pi|z^i-z^j|^2}$. Supposing for instance that $a_1<a_2$,  and fixing $z^2,\hdots,z^n\in\R^2$ all distinct with $z^2\neq 0$, one deduces that $L^*p_\infty(z)$ goes to $+\infty$ when $z^1=z^2+\lambda{z^2}^\perp$ with $\lambda\to 0^+$.
\end{arem}
The following computations about the exponential decay of the relative entropy are inspired from the proofs of Lemma 3.2 \cite{gawa} and Proposition 8
  \cite{JLLO}. 
   Since $\nabla_{z^i}. K(z^i-z^j)=0$, one can replace the
  expression
$K(z^i-z^j) \nabla_{z^i} p_t$ by $\nabla_{z^i}. \left(K(z^i-z^j)
p_t\right)$. Thus, by the
Fokker-Planck equation satisfied by $p_t$,
  Stokes' formula and remarking that
  $\int_{\R^{2n}}\nabla_z.(zp_t)\ln(p_t)=-\int_{\R^{2n}} z.\nabla_z p_t=\int_{\R^{2n}}
  2np_t$, one obtains
   \begin{align*}
      \partial_t\int_{\R^{2n}}
      p_t\ln\left(\frac{p_t}{p_\infty}\right)&=\int_{\R^{2n}}
      \partial_tp_t\ln\left(\frac{p_t}{p_\infty}\right)\\
    &=  \int_{\R^{2n}}\sum_{i=1}^n\nabla_{z^i}.\bigg(\nu\nabla_{z^i}
       p_t+\frac{z^i}{2}p_t -a\sum_{\stackrel{j=1}{j\neq i}}^nK(z^i-z^j) p_t\bigg)
        \times(\ln(p_t)-\ln(p_\infty))\\
&=-\int_{\R^{2n}}\nu\frac{|\nabla_z
  p_t|^2}{p_t}-np_t-a\sum_{j\neq
  i}K(z^i-z^j).\nabla_{z^i}p_t\\\phantom{=}&+\int_{\R^{2n}}\left(-a\sum_{j\neq i}K(z^i-z^j).\nabla_{z^i}\ln(p_\infty)+\frac{z}{2}.\nabla_z\ln(p_\infty)-\nu\Delta_z\ln(p_\infty)\right)p_t. \end{align*}
Next, dividing the stationary Fokker-Planck equation by $p_\infty$,
   one obtains
\begin{equation}
\frac{z}{2}.\nabla_z\ln(p_\infty)-a\sum_{j\neq i}K(z^i-z^j).\nabla_{z^i}\ln(p_\infty)=-\nu\left(\Delta_z\ln(p_\infty)+|\nabla_z\ln(p_\infty)|^2\right) -n.\label{statfplog}\end{equation} Using this and $\int_{\R^{2n}}\sum_{j\neq
  i}K(z^i-z^j).\nabla_{z^i}p_t=\sum_{i=1}^n\int_{\R^{2n}}\nabla_{z^i}.\left(\sum_{\stackrel{j=1}{j\neq i}}^nK(z^i-z^j)p_t\right)=0$ yields
\begin{equation}\label{dissip}
\begin{split}
\partial_t\int_{\R^{2n}} p_t\ln\left(\frac{p_t}{p_\infty}\right)&=-\nu\int_{\R^{2n}}\frac{|\nabla_z
  p_t|^2}{p_t}+
2\Delta_z\ln(p_\infty)p_t+|\nabla_z\ln(p_\infty)|^2p_t\\
&=-\nu\int_{\R^{2n}}\left|\nabla_z\ln\left(\frac{p_t}{p_\infty}\right)\right|^2p_t\\
&=-4\nu\int_{\R^{2n}}\left|\nabla_z\sqrt{\frac{p_t}{p_\infty}}\right|^2p_\infty.
   \end{split}
   \end{equation}
The Logarithmic Sobolev Inequality satisfied by the invariant density
$p_\infty$ (see \cite{gross}) ensures that
$$\int_{\R^{2n}} p_t\ln\left(\frac{p_t}{p_\infty}\right)\leq
4\nu\int_{\R^{2n}}\left|\nabla_z\sqrt{\frac{p_t}{p_\infty}}\right|^2p_\infty,$$
and one easily concludes by comparison with the ODE $\alpha'(t)=-\alpha(t)$.
   
 \begin{adem}[of Proposition \ref{convdensnorm}]  
 Since the regularized kernel $K_{\varepsilon}$  also satisfies $ K_{\varepsilon}(z^i-z^j) \cdot (z^i-z^j) =0$, 
 the arguments before Remark \ref{noninv} permit to check that 
 $$\int_{\R^{2n}} L^{\varepsilon} f(z)p_\infty(z)dz=0$$ 
 for each compactly supported smooth function $f:\R^{2n}\to \R$, where $L^{\varepsilon}$ is the generator of the process $Z^{\varepsilon}=(Z^{\varepsilon,1}_t,\hdots,Z^{\varepsilon,n}_t)_{t\geq 0}$ defined by
$$Z^{\varepsilon,i}_t 
=X^i_0+\sqrt{2\nu}B^i_t+\int_0^t\bigg(\sum_{\stackrel{j=1}{j\neq
    i}}^na_jK_\varepsilon(Z^{\varepsilon,i}_s-Z^{\varepsilon,j}_s)-\frac{Z^{\varepsilon,i}_s}{2}\bigg)ds,\;1\leq
    i\leq n.$$ 
  This implies by the criterion of Echeverria \cite{Echev}  that the density $p_\infty$ is invariant for the process $Z^{\varepsilon}$. In particular, for each  function $f$ as before,  one has 
$$\E_{\infty} (f(Z^{\varepsilon}_t))= \int_{\R^{2n}} f(z)p_\infty(z)dz,$$
where $\E_{\infty}$ denotes the expectation when the initial condition $(X^1_0,\hdots,X^n_0)$ has the density $p_{\infty}$.
Since by the proof of Lemma \ref{feller}, for each $t>0$, $Z^{\varepsilon}_t\to Z_t$ in law as $\varepsilon \to 0$  under any initial distribution  of $Z^{\varepsilon}_0=Z_0$, we deduce that 
$$\E_{\infty} (f(Z_t))= \int_{\R^{2n}} f(z)p_\infty(z)dz$$
for all such $f$ and all $t>0$  which entails the desired result.
Since the invariant density $p_{\infty}$ of $Z^\varepsilon$ has a locally integrable gradient, assumptions  
 $H1)$ , $H2)_{p_{\infty}}$,  $H3)_{p_{\infty}}$  in \cite{FJ} are satisfied.  For any initial condition $(X^1_0,\hdots,X^n_0)$, $Z^\varepsilon_t$ admits a positive density by the boundedness of $K_\varepsilon$ and Girsanov theorem. Let $p_0$ be a density on $\R^{2n}$ and $p_t^{\varepsilon}$ denote the density of $Z^\varepsilon_t$ when $(X^1_0,\hdots,X^n_0)$ has the density $p_0$. Let us first suppose that $\frac{p_0}{p_\infty}$ has polynomial growth. By Remark 2.5 \cite{FJ}, the densities $p_t^{\varepsilon}$ then satisfy assumption $H3)_{p^{\varepsilon}}$. By Theorem  2.4 and Corollary 2.6 in \cite{FJ} we obtain for all $t\geq 0$,
$$ \partial_t\int_{\R^{2n}} p_t^{\varepsilon}\ln\left(\frac{p_t^{\varepsilon}}{p_\infty}\right)\leq -4\nu\int_{\R^{2n}}\left|\nabla_z\sqrt{\frac{p_t^{\varepsilon}}{p_\infty}}\right|^2p_\infty.$$
The Logarithmic Sobolev Inequality ensures that 
 \begin{equation}\label{expdecvareps}
 \forall t\geq 0,\;\int_{\R^{2n}}
   p^{\varepsilon}_t\ln\left(\frac{p^{\varepsilon}_t}{p_\infty}\right)\leq e^{-t}\int_{\R^{2n}}
   p_0\ln\left(\frac{p_0}{p_\infty}\right)
   \end{equation}
 for all such  $p_0$.  In order to  extend the above inequality  to a general initial density $p_0$ , recall  that  the  relative entropy $H(\cdot \vert p_{\infty})$   defined on  probability measures on $\R^{2n}$ by    
\begin{equation*}
  m\mapsto  H(m|p_{\infty})=\begin{cases}
      \int_{\R^{2n}} \ln \left(\frac{p}{p_\infty}(z)\right) p(z)dz
 \ \mbox{ if
}m(dz)=p(z)dz\\
+\infty\mbox{ otherwise},
   \end{cases}
\end{equation*}
is lower semi-continuous with respect to the weak convergence.\\ 
For $k\in{\mathbb N}^*$, let $\alpha_k=\int_{\R^{2n}}(p_0(z)\wedge kp_\infty(z))dz$, $\tilde{p}_0^k=(p_0\wedge kp_\infty)/\alpha_k$, and $p_0^k=\alpha_k\tilde{p}_0^k+(1-\alpha_k)p_\infty$. By convexity of $m\mapsto H(m|p_\infty)$,
\begin{align*}
   H(p_0^k|p_\infty)\leq& \alpha_kH(\tilde{p}_0^k|p_\infty)+(1-\alpha_k)\times 0\\
=&-\alpha_k\ln(\alpha_k)+\int_{\R^{2n}}p_0(z)\ln\left(\frac{p_0}{p_\infty}(z)\right)1_{\{p_0(z)\leq p_\infty(z)\}}dz\\&+\int_{\R^{2n}}(p_0(z)\wedge kp_\infty(z))\ln\left(\frac{p_0}{p_\infty}(z)\wedge k\right)1_{\{p_0(z)>p_\infty(z)\}}dz\end{align*}
where, by a slight abuse of notations, for a probability density $p$ on $\R^{2n}$, $H(p|p_\infty)$ denotes the relative entropy of $p(z)dz$ w.r.t. $p_\infty(z)dz$. As $k\to\infty$, the third term in the r.h.s. converges to $\int_{\R^{2n}}p_0(z)\ln\left(\frac{p_0}{p_\infty}(z)\right)1_{\{p_0(z)>p_\infty(z)\}}dz$ by monotone convergence and $\alpha_k\to 1$ by Lebesgue's theorem. With the lower semi-continuity of the relative entropy, this ensures that $H(p^k_0|p_\infty)\rightarrow H(p_0|p_\infty)$.\\
Denoting by $p^{k,\varepsilon}_t$ the density of $Z^{\varepsilon}_t$ when $(X^1_0,\hdots,X^n_0)$ admits the density $p^k_0$, one has by the semi-group property $$\int_{\R^{2n}}|p^\varepsilon_t(z)-p^{k,\varepsilon}_t(z)|dz\leq \int_{\R^{2n}}|p_0(z)-p^k_0(z)|dz=2(1-\alpha_k)$$ so that 
$H(p^{\varepsilon}_t|p_\infty)\leq \liminf_{k\to\infty}H(p^{k,\varepsilon}_t|p_\infty)$. Hence \eqref{expdecvareps} holds for general $p_0$. Finally,  using the  convergence in law $Z^{\varepsilon}_t\to Z_t$  as $\varepsilon \to 0$ following from  the proof of Lemma \ref{feller}  and once again the lower semicontinuity of the relative entropy, we conclude  that the required inequality holds for $(Z_t)$ .
 
 \end{adem}

\section{Two particles case : $n=2$}\label{toy}
This section is devoted to the simple case where only two vortices interact :
\begin{equation}
   \begin{cases}X^1_t=X^1_0+\sqrt{2\nu} W^1_t+\int_0^t
     a_2K(X^1_s-X^2_s)ds\\

      X^2_t=X^2_0+\sqrt{2\nu} W^2_t+\int_0^t
     a_1K(X^2_s-X^1_s)ds
   \end{cases}.\label{syst2part}
\end{equation}
Existence and uniqueness for this system follows from \cite{Tak} when $a_1a_2>0$ and from \cite{O2} for general vorticities.
After the change of variables $Z^i_t=e^{-t/2}X^i_{e^t-1}$ for $i\in\{1,2\}$, one obtains
the following dynamics :
\begin{equation}
   \begin{cases}Z^1_t=X^1_0+\sqrt{2\nu} B^1_t+\int_0^t
     a_2K(Z^1_s-Z^2_s)ds-\frac{1}{2}\int_0^t Z^1_sds\\

      Z^2_t=X^2_0+\sqrt{2\nu} B^2_t+\int_0^t
     a_1K(Z^2_s-Z^1_s)ds-\frac{1}{2}\int_0^t Z^2_sds
   \end{cases}.\label{syst2partchgt}
\end{equation}
The difference $Z_t=Z^1_t-Z^2_t$ solves 
\begin{equation}
   Z_t=X^1_0-X^2_0+2\sqrt{\nu} B_t+a\int_0^t
     K(Z_s)ds-\frac{1}{2}\int_0^t Z_sds
\label{eqz}\end{equation}
where $a\stackrel{\rm def}{=}a_1+a_2$ and 
$B_t\stackrel{\rm def}{=}\frac{B^1_t-B^2_t}{\sqrt{2}}$ is a
two-dimensional Brownian motion independent from the two dimensional
Brownian motion $B^{\perp}_t\stackrel{\rm
  def}{=}\frac{B^1_t+B^2_t}{\sqrt{2}}$.
Remarking that
\begin{equation}
   \begin{cases}
      Z^1_t=e^{-\frac{t}{2}}X^1_0+\sqrt{\nu}\int_0^te^{\frac{s-t}{2}}dB_s+\sqrt{\nu}\int_0^te^{\frac{s-t}{2}}dB^\perp_s+a_2\int_0^te^{\frac{s-t}{2}}K(Z_s)ds\\
   Z^2_t=e^{-\frac{t}{2}}X^2_0-\sqrt{\nu}\int_0^te^{\frac{s-t}{2}}dB_s+\sqrt{\nu}\int_0^te^{\frac{s-t}{2}}dB^\perp_s-a_1\int_0^te^{\frac{s-t}{2}}K(Z_s)ds\end{cases},\label{z1z2}
\end{equation}
we see that to understand the long-time behaviour of $(Z^1_t,Z^2_t)$, it
is enough to study the long-time behaviour of the triplet
$(\int_0^te^{\frac{s-t}{2}}dB_s,\int_0^te^{\frac{s-t}{2}}K(Z_s)ds,\int_0^te^{\frac{s-t}{2}}dB^\perp_s)$.
The last coordinate is independent from the two first and converges in
law to a two-dimensional standard normal random variable. So we only
need to study the long-time behaviour of the couple
$(\int_0^te^{\frac{s-t}{2}}dB_s,\int_0^te^{\frac{s-t}{2}}K(Z_s)ds)$. In
fact,
\begin{equation}
   \forall t\geq
0,\;Z_t=(X^1_0-X^2_0)e^{-\frac{t}{2}}+2\sqrt{\nu}\int_0^te^{\frac{s-t}{2}}dB_s+a\int_0^te^{\frac{s-t}{2}}K(Z_s)ds\label{eqz}
\end{equation}
and we are going to study the limit behaviour of $$\mu_t\stackrel{\rm def}{=}{\mbox Law}\left(Z_t,\int_0^te^{\frac{s-t}{2}}dB_s,\int_0^te^{\frac{s-t}{2}}K(Z_s)ds\right).$$ Let $q_\infty:z\in\R^2\mapsto\frac{1}{8\nu\pi}e^{-\frac{|z|^2}{8\nu}}$
be the density of the law ${\cal N}_2(0,4\nu I_2)$ and $\bar{Z}_0$ denote a $\R^2$-valued random variable distributed according to this law independent from the Brownian motions $B^1$ and $B^2$. The stochastic
differential equation
\begin{equation}
   \bar{Z}_t={\bar Z}_0+2\sqrt{\nu} B_t-a\int_0^t
     K(\bar{Z}_s)ds-\frac{1}{2}\int_0^t \bar{Z}_sds\label{defzbar}
\end{equation}
admits a unique solution. Indeed, for all $\varepsilon>0$, $K$ is globally Lipschitz
continuous and bounded on $\{z\in\R^2:|z|\geq \varepsilon\}$ and
$\bar{R}_t\stackrel{\rm def}{=}|\bar{Z}_t|^2$ solves the Cox-Ingersoll-Ross
stochastic differential equation
$$d\bar{R}_t=4\sqrt{\nu \bar{R}_t}d\beta_t+(8\nu-\bar{R}_t)dt$$
with $\beta_t=\int_0^t\frac{\bar{Z}_s.dB_s}{\sqrt{\bar{R}_s}}$ a standard Brownian motion and does not
vanish since $8\nu=\frac{(4\sqrt{\nu})^2}{2}$ (see Proposition 6.2.3 \cite{lamblap}).
By an easy adaptation of the proof of the first statement in Proposition \ref{convdensnorm}, one checks that the density
$q_\infty$ is stationary for both the SDEs satisfied by $Z_t$ and
$\bar{Z}_t$. As a consequence, for all $s\geq 0$, $\bar{Z}_s$ admits the
density $q_\infty$ and $\E(|\bar{Z}_s|)=\E(|\bar{Z}_0|)<+\infty$,
$\E(|K(\bar{Z}_s)|)=\frac{1}{2\pi}\E(|\bar{Z}_0|^{-1})<+\infty$. One
deduces that the random vector $$\left(\bar{Z}_0,\int_0^\infty
e^{-\frac{s}{2}}\left(\bar{Z}_s/\sqrt{4\nu}\,ds-dB_s\right),\int_0^\infty
e^{-\frac{s}{2}}K(\bar{Z}_s)ds\right)$$ is well defined. Let $\mu_\infty$
denote its distribution.
\begin{aprop}
   If $X^1_0-X^2_0$ admits a density $q_0$ with respect to the
   Lebesgue measure on $\R^2$ such that
   $\int_{\R^2}q_0\ln\left(\frac{q_0}{q_\infty}\right)<+\infty$, then
   for all $t\geq 0$, the density $q_t$ of $Z_t$ is such that
\begin{equation}
   \int_{\R^2}q_t\ln\left(\frac{q_t}{q_\infty}\right)\leq e^{-t}\int_{\R^2}q_0\ln\left(\frac{q_0}{q_\infty}\right).\label{convz}
\end{equation}
Moreover, for all $\alpha\in[1,2)$, $W_\alpha(\mu_t,\mu_\infty)$ converges to
   $0$ exponentially fast as $t\to +\infty$. If, in addition,
   $\E(|X^1_0-X^2_0|^\rho)<+\infty$ for some $\rho>2$, then this
   exponential convergence still holds for all $\alpha\in[1,\rho)$.
\label{convloidif}\end{aprop}
\begin{arem}\label{transp}
   By the transport inequality satisfied by the Gaussian law (see \cite{Tal}), $W_2(q_0,q_\infty)\leq
   2\sqrt{\nu\int_{\R^2}q_0\ln\left(\frac{q_0}{q_\infty}\right)}$
   (where by a slight abuse of notations $q_\infty$ and $q_0$ stand for
   the measures $q_\infty(z)dz$ and $q_0(z)dz$). Hence, the finiteness
   of the relative entropy
   $\int_{\R^2}q_0\ln\left(\frac{q_0}{q_\infty}\right)$ implies the
   finiteness of $\E(|X^1_0-X^2_0|^2)$.
\end{arem}
This proposition, the proof of which is postponed, is the main step in the derivation of the long-time
behaviour of $(Z^1_t,Z^2_t)$.
\begin{thm}If $X^1_0-X^2_0$ admits a density $q_0$ such that
   $\int_{\R^2}q_0\ln\left(\frac{q_0}{q_\infty}\right)<+\infty$, then the $2\times 2$ matrix $(Z^1_t,Z^2_t)$ (with first column equal to $Z^1_t$ and second column equal to $Z^2_t$) converges in law to 
\begin{align*}
  \frac{\bar{Z}_0}{2}(1,-1)+\int_0^{+\infty}e^{-\frac{s}{2}}\left(\sqrt{\nu}dB^\perp_s+\frac{(a_2-a_1)}{2}K(\bar{Z}_s)ds\right)(1,1).
\end{align*}
Let $\mu^{1,2}_t$ and $\mu^{1,2}_\infty$ respectively denote the law of
$(Z^1_t,Z^2_t)$ and of the above limit. For $\alpha\geq 1$ such that $\E(|X^1_0|^\alpha+|X^2_0|^\alpha)<+\infty$, $W_\alpha(\mu^{1,2}_t,\mu^{1,2}_\infty)$ converges
to $0$ exponentially fast when $\alpha<2$ and under the additional condition $\E(|X^1_0-X^2_0|^\rho)<+\infty$ for some $\rho>\alpha$ when $\alpha\geq 2$.\\% , then $W_\alpha(\mu^{1,2}_t,\mu^{1,2}_\infty)$ converges
% to $0$ exponentially fast for all
% $\alpha\in[1,\rho)$ and, in case $\rho<2$, also for $\alpha=\rho$.\\
Last, unless $a_2=a_1$, the limiting distribution $\mu^{1,2}_\infty$ is not Gaussian.
\label{convloi12}\end{thm}
\begin{arem}\begin{itemize}
   \item The linear combination $a_1Z^1_t+a_2Z^2_t$ is an
  Ornstein-Uhlenbeck process and converges in law to some Gaussian
  limit. The difference $Z_t=Z^1_t-Z^2_t$, which is linearly independent
  of $a_1Z^1_t+a_2Z^2_t$ as soon as $a_1+a_2\neq 0$, also converges to
  some Gaussian limit. Nethertheless, unless $a_1=a_2$, the limit distribution of
$(Z^1_t,Z^2_t)$ is not Gaussian.
\item When $a_1=a_2$, $(Z^1_t,Z^2_t)$ converge in distribution to 
$\frac{\bar{Z}_0}{2}(1,1)+\int_0^{+\infty}e^{-\frac{s}{2}}\sqrt{\nu}dB^\perp_s(1,-1)$ with $\frac{\bar{Z_0}}{2}$ and $\int_0^{+\infty}e^{-\frac{s}{2}}\sqrt{\nu}dB^\perp_s$ both distributed according to ${\cal N}_2(0,\nu I_2)$ and independent. One thus recovers the limit distribution obtained in 
Proposition \ref{convdensnorm}.
\end{itemize}\end{arem}

Let us first check that the Theorem follows from Proposition
\ref{convloidif} before proving this Proposition.
\begin{adem}[of Theorem \ref{convloi12}]
Let $\alpha\geq 1$. The random vector $\int_0^te^{\frac{s-t}{2}}dB^\perp_s$ distributed
according to ${\cal N}_2(0,(1-e^{-t})I_2)$ converges in law to $G\sim{\cal
  N}_2(0,I_2)$. Moreover, $W_\alpha({\cal N}_2(0,(1-e^{-t})I_2),{\cal
  N}_2(0,I_2))\leq
e^{-\frac{t}{2}}\left(\E(|G|^\alpha)\right)^{1/\alpha}$ converges
exponentially fast to $0$. Since for $x\in\R^6$ and $y\in\R^2$, $(|x|^2+|y|^2)^{\alpha/2}\leq 2^{(\frac{\alpha}{2}-1)^+}(|x|^\alpha+|y|^\alpha)$, one gets $$W_\alpha\left(\mu_t\otimes{\cal
  N}_2(0,(1-e^{-t})I_2),\mu_\infty\otimes {\cal
  N}_2(0,I_2)\right)\leq 2^{(\frac{1}{2}-\frac{1}{\alpha})^+}\left(W_\alpha(\mu_t,\mu_\infty)+W_\alpha({\cal N}_2(0,(1-e^{-t})I_2),{\cal
  N}_2(0,I_2))\right)$$
since, for $\rho_1$ and $\rho_2$ optimal couplings respectively between $\mu_t$ and $\mu_\infty$ and between ${\cal
  N}_2(0,(1-e^{-t})I_2)$ and ${\cal
  N}_2(0,I_2)$, the probability measure $\rho(dxdy,d\tilde{x}d\tilde{y})=\rho_1(dx,d\tilde{x})\rho_2(dy,d\tilde{y})$ is a coupling between $\mu_t\otimes{\cal
  N}_2(0,(1-e^{-t})I_2)$ and $\mu_\infty\otimes{\cal
  N}_2(0,I_2)$.
Hence, by Proposition \ref{convloidif}, the Wasserstein distance between $$\mu_t\otimes{\cal
  N}_2(0,(1-e^{-t})I_2)={\cal L}aw\left(Z_t,\int_0^te^{\frac{s-t}{2}}dB_s,\int_0^te^{\frac{s-t}{2}}K(Z_s)ds,\int_0^te^{\frac{s-t}{2}}dB^\perp_s\right)$$
and $$\mu_\infty\otimes{\cal
  N}_2(0,I_2)={\cal L}aw\left(\bar{Z}_0,\int_0^\infty
e^{-\frac{s}{2}}(\bar{Z}_s/\sqrt{4\nu}\,ds-dB_s),\int_0^\infty
e^{-\frac{s}{2}}K(\bar{Z}_s)ds,\int_0^\infty
e^{-\frac{s}{2}}dB^\perp_s \right)$$ converges exponentially fast to $0$  when $\alpha<2$ and under the condition $\E(|X^1_0-X^2_0|^\rho)<+\infty$ for some $\rho>\alpha$ when $\alpha\geq 2$.
The first statement then follows from \eqref{z1z2} and Slutsky's theorem
: $(Z^1_t,Z^2_t)$ writes as the sum of $e^{-\frac{t}{2}}(X^1_0,X^2_0)$
which converges a.s. to zero and the vector 
$$\left(\sqrt{\nu}\int_0^t e^{\frac{s-t}{2}}(dB_s+dB^\perp_s)+a_2\int_0^t e^{\frac{s-t}{2}}K(Z_s)ds,\sqrt{\nu}\int_0^t e^{\frac{s-t}{2}}(dB^\perp_s-dB_s)-a_1\int_0^t e^{\frac{s-t}{2}}K(Z_s)ds\right)$$
with law $\tilde{\mu}^{1,2}_t$ which converges weakly to
$$\int_0^{+\infty}e^{-\frac{s}{2}}\left(\sqrt{\nu}dB^\perp_s+\frac{(a_2-a_1)}{2}K(\bar{Z}_s)ds\right)(1,1)+\int_0^{+\infty}e^{-\frac{s}{2}}\left(-\sqrt{\nu}dB_s+\frac{1}{2}\bar{Z}_sds+\frac{a}{2}K(\bar{Z}_s)ds\right)(1,-1),$$
where
$e^{-\frac{s}{2}}\left(-\sqrt{\nu}dB_s+\frac{1}{2}\bar{Z}_sds+\frac{a}{2}K(\bar{Z}_s)ds\right)=-\frac{1}{2}d(e^{-\frac{s}{2}}\bar{Z}_s)$
by \eqref{defzbar}.\\
The second statement is a consequence of the triangle inequality 
$$W_\alpha(\mu^{1,2}_t,\mu^{1,2}_\infty)\leq
W_\alpha(\mu^{1,2}_t,\tilde{\mu}^{1,2}_t)+W_\alpha(\tilde{\mu}^{1,2}_t,\mu^{1,2}_\infty)$$
where by \eqref{z1z2} and the above definition of $\tilde{\mu}^{1,2}_t$, the first term of the right-hand-side is not greater than
$e^{-\frac{t}{2}}\left(\E(|(X^{1}_0,X^{2}_0)|^\alpha\right)^{1/\alpha}$.
Since $\tilde{\mu}^{1,2}_t$ (resp. $\tilde{\mu}^{1,2}_\infty$) is the image of $\mu_t\otimes{\cal
  N}_2(0,(1-e^{-t})I_2)$ (resp. $\mu_\infty\otimes {\cal
  N}_2(0,I_2)$) by the mapping $x\in\R^8\mapsto \left(\begin{array}{cccccccc}0 & 0 &\sqrt{\nu}& 0 & a_2 &0 &\sqrt{\nu}& 0\\0& 0& 0 &\sqrt{\nu}& 0&a_2 & 0&\sqrt{\nu}\\0& 0 &-\sqrt{\nu}& 0& -a_1 & 0&\sqrt{\nu}& 0\\0& 0 & 0 &-\sqrt{\nu}& 0& -a_1 & 0&\sqrt{\nu}
\end{array}\right)x$, the second term is not greater than $W_\alpha\left(\mu_t\otimes{\cal
  N}_2(0,(1-e^{-t})I_2),\mu_\infty\otimes {\cal
  N}_2(0,I_2)\right)$ multiplied by the operator norm of the above matrix.

Let us now suppose that $a_1\neq a_2$ and that $\mu_\infty^{1,2}$ is Gaussian and obtain some contradiction. Then, as can be seen using the characteristic function, the limit $(\bar{Z}_0,\int_0^{+\infty}e^{-\frac{s}{2}}K(\bar{Z}_s)ds)$ of $(Z_t,\int_0^te^{\frac{s-t}{2}}K(Z_s)ds)$ is Gaussian by independence of $(\bar{Z}_t)_{t\geq 0}$ and $\int_0^\infty e^{-\frac{s}{2}}dB^\perp_s$.\\Let $\zeta_t=2\sqrt{\nu}\int_0^te^{\frac{s-t}{2}}dB_s\sim{\cal N}_2(0,4\nu(1-e^{-t}) I_2)$ which solves $d\zeta_t=2\sqrt{\nu}dB_t-\frac{\zeta_t}{2}dt$. Since by \eqref{eqz}, $$\zeta_t=Z_t-e^{-\frac{t}{2}}Z_0-a\int_0^te^{\frac{s-t}{2}}K(Z_s)ds,$$ $(Z_t,\zeta_t)$ converges in law to some Gaussian random vector $(Z_\infty,\zeta_\infty)$ with $Z_\infty$ and $\zeta_\infty$ both distributed according to ${\cal N}_2(0,4\nu I_2)$. Let $\beta$ and $\gamma$ be the respective correlations of $\zeta_{1,\infty}$ with $Z_{1,\infty}$ and $Z_{2,\infty}$ and $$L=2\nu((\partial_{z_1}+\partial_{\zeta_1})(\partial_{z_1}+\partial_{\zeta_1})+(\partial_{z_2}+\partial_{\zeta_2})(\partial_{z_2}+\partial_{\zeta_2}))+(aK(z)-\frac{z}{2}).\nabla_z-\frac{\zeta}{2}.\nabla_\zeta$$ where $z=(z_1,z_2)$ and $\zeta=(\zeta_1,\zeta_2)$ be the infinitesimal generator of the process $(Z_t,\zeta_t)$. When both $f:\R^4\to\R$ and $Lf$ are continuous and bounded functions, then $\E(f(Z_t,\zeta_t))$ and $\E(Lf(Z_t,\zeta_t))$ converge respectively to $\E(f(Z_\infty,\zeta_\infty))$ and $\E(Lf(Z_\infty,\zeta_\infty))$ as $t\to\infty$. If moreover $\nabla f$ is bounded, taking expectations in Itô's formula, one obtains $$\frac{d}{dt}\E(f(Z_t,\zeta_t))=\E(Lf(Z_t,\zeta_t)),$$ so that the previous convergences ensure that $\E(Lf(Z_\infty,\zeta_\infty))=0$.
Assume for a while that this centering property may be extended to the choices $f(z,\zeta)$ respectively equal to $\zeta_1z_1$, $\zeta_1z_2$ and $\zeta_1^2z_1^2$, for which $Lf$ is equal to $4\nu-z_1\zeta_1-a\frac{z_2\zeta_1}{2\pi|z|^2}$, $-\zeta_1z_2+a\frac{\zeta_1z_1}{2\pi|z|^2}$ and $4\nu(\zeta_1^2+z_1^2+4\zeta_1z_1)-2z_1^2\zeta_1^2-\frac{2az_1z_2\zeta_1^2}{2\pi|z|^2}$ and $\E(Lf(Z_\infty,\zeta_\infty))$ to $4\nu(1-\beta)-\frac{a\gamma}{4\pi}$, $-4\nu\gamma +\frac{a\beta}{4\pi}$ and $(4\nu)^2(2+4\beta)-2(4\nu)^2(2\beta^2+1)-\frac{2\nu}{\pi} a\beta\gamma$. One obtains the system \begin{equation}
   \left\{\begin{array}{l}4\nu(1-\beta)-\frac{a\gamma}{4\pi}=0\\-4\nu\gamma +\frac{a\beta}{4\pi}=0\\\beta\left(16\nu(1-\beta)-\frac{a\gamma}{2\pi}\right)=0
\end{array}\right.\Leftrightarrow\left\{\begin{array}{l}a\beta\gamma=0\\
-4\nu\gamma +\frac{a\beta}{4\pi}=0\\
4\nu(1-\beta)-\frac{a\gamma}{4\pi}=0\end{array}\right.\label{systincomp}
\end{equation}
which has no solution $(\beta,\gamma)$ when $a\neq 0$. This provides the required contradiction in the case $a\neq 0$, once we are able to extend the centering property to the three above choices of $f$. To deal with the singularity of the terms like $\frac{z_2\zeta_1}{2\pi|z|^2}$, for $\varepsilon\in(0,1)$, we introduce a non-decreasing odd $C^2$ cutoff function $\varphi_\varepsilon:\R\to\R$ such that $$\varphi_\varepsilon(x)=\begin{cases}0\mbox{ if }|x|\leq \frac{\varepsilon}{2}\\
x\mbox{ if }|x|\in[\varepsilon,\frac{1}{\varepsilon}]\\
\frac{3x}{2\varepsilon|x|}\mbox{ if }x\geq\frac{2}{\varepsilon}
\end{cases},\;\|\varphi'_\varepsilon\|_\infty\leq 2\mbox{ and }\sup_{\varepsilon<1,|x|\geq \varepsilon}|\varphi''_\varepsilon(x)|+\varepsilon\sup_{\varepsilon<1,|x|\leq \varepsilon}|\varphi''_\varepsilon(x)|<\infty.$$
The centering property $\E(Lf_\varepsilon(Z_\infty,\zeta_\infty))=0$ holds for $f_\varepsilon$ respectively equal to $\zeta_1\varphi_\varepsilon(z_1)$, $\zeta_1\varphi_\varepsilon(z_2)$ and $\zeta_1^2\varphi^2_\varepsilon(z_1)$ by replacing the boundedness by the uniform integrability derived from the bound $$\forall n\in\N,\,\sup_{t\geq 0}\E(|\zeta_t|^n)=\E(|\zeta_\infty|^n)<+\infty.$$ For $f(z,\zeta)=\zeta_1z_1$ and $f_\varepsilon(z,\zeta)=\zeta_1\varphi_\varepsilon(z_1)$, one has
$$L(f_\varepsilon-f)(Z_\infty,\zeta_\infty)=2\nu\zeta_{1,\infty}\varphi_\varepsilon''(Z_{1,\infty})+(1-\varphi_\varepsilon'(Z_{1,\infty}))\times \left(-4\nu+\frac{Z_{1,\infty}\zeta_{1,\infty}}{2}+a\frac{Z_{2,\infty}\zeta_{1,\infty}}{2\pi|Z_\infty|^2}\right).$$
The expectation of the second term in the right-hand-side tends to $0$ as $\varepsilon\to 0$ by Lebesgue's theorem. On the other hand, using the controls on $\varphi_\varepsilon''$ for the inequality, one obtains
\begin{align*}
   |\E(\zeta_{1,\infty}\varphi_\varepsilon''(Z_{1,\infty}))|&=|\E(\E(\zeta_{1,\infty}\varphi_\varepsilon''(Z_{1,\infty}))|Z_{1,\infty})|=|\beta \E(Z_{1,\infty}\varphi_\varepsilon''(Z_{1,\infty}))|\\
&\leq C\left(\E\left(\frac{|Z_{1,\infty}|}{\varepsilon}1_{\{\frac{\varepsilon}{2}\leq |Z_{1,\infty}|\leq \varepsilon\}}\right)+\E\left(|Z_{1,\infty}|1_{\{|Z_{1,\infty}|\geq\frac{1}{\varepsilon}\}}\right)\right)\stackrel{\varepsilon\to 0}{\longrightarrow}0,
\end{align*}
so that $\E(Lf(Z_\infty,\zeta_\infty))=0$. By similar arguments, the centering is still true for $f$ equal to $\zeta_1z_2$ and $\zeta_1^2z_1^2$

When $a=0$, $(\beta,\gamma)=(1,0)$ solves the system \eqref{systincomp}, which is not surprising since $\forall t,\;Z_t=\zeta_t$. We then work with $\xi_t=\int_0^te^{\frac{s-t}{2}}K(Z_s)ds$. As $t\to\infty$,  $(Z_t,\xi_t)$ converges in distribution to some Gaussian random vector $(Z_\infty,\xi_\infty)$ with $Z_\infty$ distributed according to ${\cal N}_2(0,4\nu I_2)$ and $\xi_\infty$ centered since $\int_0^\infty e^{-\frac{s}{2}}K(\bar{Z}_s)ds$ is centered. The infinitesimal generator of $(Z_t,\xi_t)$ is $L=2\nu\Delta_z-\frac{z}{2}.\nabla_z+(K(z)-\frac{\xi}{2}).\nabla_\xi$. For $f(z,\xi)$ respectively equal to $-\xi_1z_2$ and $\xi^2_1 z^2_2$, $Lf$ is equal to $\xi_1z_2+\frac{z_2^2}{2\pi|z|^2}$ and $4\nu\xi_1^2-2\xi^2_1 z^2_2-\frac{\xi_1z^3_2}{\pi|z|^2}$ and the equality $\E(Lf(Z_\infty,\xi_\infty))=0$ yields \begin{equation*}
   \begin{cases}
      {\rm Cov}(Z_{2,\infty},\xi_{1,\infty})=-\frac{1}{4\pi}\\4\nu{\rm Var}(\xi_{1,\infty})-2\left(4\nu{\rm Var}(\xi_{1,\infty})+2{\rm Cov}^2(Z_{2,\infty},\xi_{1,\infty})\right)-\frac{3}{4\pi}{\rm Cov}(Z_{2,\infty},\xi_{1,\infty})=0
   \end{cases}.
\end{equation*} These equalities imply ${\rm Var}(\xi_{1,\infty})=-\frac{1}{64\nu\pi^2}$ which is the desired contradiction.
To justify the equality $\E(Lf(Z_\infty,\xi_\infty))=0$ for the above choices of $f$ one first construct approximations $f_{\varepsilon,\eta}$ for  $\varepsilon,\eta\in (0,1)$ by replacing the factors $z_1$, $z_2$, $\xi_1$, $\xi_2$ respectively by $\varphi_\varepsilon(z_1)$, $\varphi_\varepsilon(z_2)$, $\varphi_\eta(\xi_1)$ and $\varphi_\eta(\xi_2)$. For $\varepsilon\in (0,1)$, $f_\varepsilon$ is obtained similarly by only replacing the factors $z_1$ and $z_2$ in $f$. Then $\E(Lf_{\varepsilon,\eta}(Z_\infty,\xi_\infty))=0$. Since no second order derivative of $\varphi_\eta$ appears in $Lf_{\varepsilon,\eta}$, one may apply Lebesgue's theorem to take the limit $\eta\to 0$ and obtain $\E(Lf_\varepsilon(Z_\infty,\xi_\infty))=0$. Then one concludes as previously by taking the limit $\varepsilon\to 0$.
  
\end{adem}

\begin{adem}[of Proposition \ref{convloidif}]
The first statement is obtained by an easy adaptation of the proof of
Proposition \ref{convdensnorm}. Let us deal with the second statement.
Let $t>0$. Csiszar-Kullback inequality writes $\|q_{t/2}-q_\infty\|_{1}\leq
\sqrt{2\int_{\R^2}q_{t/2}\ln\left(\frac{q_{t/2}}{q_\infty}\right)}$. By \eqref{convz}, the left-hand-side is smaller than $Ce^{-\frac{t}{4}}$ with
$C$ not depending on $t$. Let $\rho_{t/2}=\frac{q_{t/2}\wedge q_\infty}{q_{t/2}}$ and
$(U,\zeta_{t/2})$ be a couple of independent random variables
independent from $(B^1,B^2,X^1_0,X^2_0)$ with $U$ uniformly distributed
on $[0,1]$ and $\zeta_{t/2}$ distributed according to the density
$\frac{(q_\infty-q_{t/2})^+}{\int_{\R^2}(q_\infty-q_{t/2})^+}$ (when
$\int_{\R^2}(q_\infty-q_{t/2})^+=0$, $q_{t/2}=q_\infty$ and
$\zeta_{t/2}$ is not needed in what follows). The random variable
$Z^{t/2}_{t/2}=1_{\{U\leq
  \rho_{t/2}(Z_{t/2})\}}Z_{t/2}+1_{\{U>\rho_{t/2}(Z_{t/2})\}}\zeta_{t/2}$ admits the density $q_\infty$ and is such that $\P(Z_{t/2}\neq Z^{t/2}_{t/2})=\frac{1}{2}\|q_{t/2}-q_\infty\|_{1}$.
Let $(Z^{t/2}_s)_{s\geq t/2}$ be the unique solution of the SDE
$$Z^{t/2}_s=Z^{t/2}_{t/2}+2\sqrt{\nu}(B_s-B_{t/2})+\int_{t/2}^s
\left(aK(Z^{t/2}_r)-\frac{1}{2}Z^{t/2}_r\right)dr$$
and $\tilde{\mu}_t$ denote the law of
$(Z^{t/2}_t,\int_{t/2}^te^{\frac{s-t}{2}}dB_s,\int_{t/2}^te^{\frac{s-t}{2}}K(Z^{t/2}_s)ds)$.
By trajectorial uniqueness, $(Z_s)_{s\geq t/2}$ and $(Z^{t/2}_s)_{s\geq
  t/2}$ coincide on $\{U\leq \rho_{t/2}(Z_{t/2})\}$. Let $A_{t/2}=\{U>\rho_{t/2}(Z_{t/2})\}$.
By the triangle inequality, for $\alpha\geq 1$,
\begin{equation}
   W_\alpha(\mu_t,\mu_\infty)\leq
W_\alpha(\mu_t,\tilde{\mu}_t)+W_\alpha(\tilde{\mu}_t,\mu_\infty).\label{triangle}
\end{equation}
For $\gamma>\alpha$, one has, using H\"older's inequality for the second step,
\begin{align*}
   &W_\alpha^\alpha(\mu_t,\tilde{\mu}_t)\leq
   \E\bigg(\bigg(1_{A_{t/2}}|Z_t-Z^{t/2}_t|^2+\bigg|\int_0^{t/2}e^{\frac{s-t}{2}}dB_s\bigg|^2\\&\phantom{W_\alpha^\alpha(\mu_t,\tilde{\mu}_t)\leq
   \E\bigg(\bigg(}+\bigg|\int_0^{t/2}e^{\frac{s-t}{2}}K(Z_s)ds+1_{A_{t/2}}\int_{t/2}^te^{\frac{s-t}{2}}(K(Z_s)-K(Z^{t/2}_s))ds\bigg|^2\bigg)^{\alpha/2}\bigg)\\
&\leq
3^{(\frac{\alpha}{2}-1)^+}\bigg[2^{\alpha-1}\P(A_{t/2})^{\frac{\gamma-\alpha}{\gamma}}(\E^{\frac{\alpha}{\gamma}}|Z_t|^\gamma+\E^{\frac{\alpha}{\gamma}}|Z^{t/2}_t|^\gamma)+\E\bigg|\int_0^{t/2}e^{\frac{s-t}{2}}dB_s\bigg|^\alpha\\
&+3^{\alpha-1}\bigg\{\E\bigg|\int_0^{t/2}e^{\frac{s-t}{2}}K(Z_s)ds\bigg|^\alpha+\P(A_{t/2})^{\frac{\gamma-\alpha}{\gamma}}\bigg(\E^{\frac{\alpha}{\gamma}}\bigg|\int_{t/2}^te^{\frac{s-t}{2}}K(Z_s)ds\bigg|^\gamma+\E^{\frac{\alpha}{\gamma}}\bigg|\int_{t/2}^te^{\frac{s-t}{2}}K(Z^{t/2}_s)ds\bigg|^\gamma\bigg)\bigg\}\bigg]
\end{align*}
The term
$\E\bigg|\int_0^{t/2}e^{\frac{s-t}{2}}dB_s\bigg|^\alpha$ is equal to $(e^{-t/2}-e^{-t})^{\alpha/2}\E|G|^\alpha$
where $G\sim{\cal N}_2(0,I_2)$ and converges exponentially fast to
$0$. So does $\P(A_{t/2})^{\frac{\gamma-\alpha}{\gamma}}=\left(\frac{1}{2}\|q_{t/2}-q_\infty\|_{1}\right)^{\frac{\gamma-\alpha}{\gamma}}$ according to the
beginning of the proof. According to Lemma \ref{contmominv} below and since $\forall z\in\R^2$, $|K(z)|\leq\frac{2\pi}{|z|}$,  the terms
involving the Biot and Savart
kernel $K$ also converge to $0$ exponentially fast. The expectation $\E|Z^{t/2}_t|^\gamma$ does not
depend on $t$ and is finite. When
$\int_{\R^2}q_0\ln\left(\frac{q_0}{q_\infty}\right)<+\infty$, choosing
$\gamma=2$, we deduce from Remark \ref{transp} and Lemma \ref{contmomz}
below that
$W_\alpha(\mu_t,\tilde{\mu}_t)$ converges to $0$ exponentially fast for
$\alpha\in[1,2)$.  When moreover $\E|X^1_0-X^2_0|^\rho<+\infty$ for
some $\rho>2$, choosing $\gamma=\rho$, we obtain that this exponential
convergence holds for
$\alpha\in[1,\rho)$.

By Lemma \ref{timerev} below, which is based on a time-reversal argument,
$\tilde{\mu}_t$ is the law of
$(\bar{Z}_0,\int_0^{t/2}e^{-\frac{s}{2}}(\bar{Z}_s/\sqrt{4\nu}\,ds-dB_s),\int_0^{t/2}e^{-\frac{s}{2}}K(\bar{Z}_s)ds)$.
Since $\mu_\infty$ is the law of $(\bar{Z}_0,\int_0^{\infty}e^{-\frac{s}{2}}(\bar{Z}_s/\sqrt{4\nu}\,ds-dB_s),\int_0^{\infty}e^{-\frac{s}{2}}K(\bar{Z}_s)ds)$, one deduces that
$$W_\alpha^\alpha(\tilde{\mu}_t,\mu_\infty)\leq 2^{(\frac{\alpha}{2}-1)^+}\left(\frac{1}{2\nu^{\alpha/2}}\E\bigg|\int_{t/2}^{+\infty}e^{-\frac{s}{2}}\bar{Z}_sds\bigg|^\alpha+2^{\alpha-1}\E\bigg|\int_{t/2}^{+\infty}e^{-\frac{s}{2}}dB_s\bigg|^\alpha+\E\bigg|\int_{t/2}^{+\infty}e^{-\frac{s}{2}}K(\bar{Z}_s)\bigg|^\alpha\right).$$
The term
$\E\bigg|\int_{t/2}^{+\infty}e^{-\frac{s}{2}}dB_s\bigg|^\alpha$ is equal to
$e^{-\frac{\alpha t}{4}}\E|G|^\alpha$ 
where $G\sim{\cal N}_2(0,I_2)$ and converges exponentially fast to
$0$. Moreover, by Hölder's inequality 
$$\E\bigg|\int_{t/2}^{+\infty}e^{-\frac{s}{2}}\bar{Z}_sds\bigg|^\alpha\leq
2^{\alpha-1}e^{-\frac{(\alpha-1)t}{4}}\int_{t/2}^{+\infty}e^{-\frac{s}{2}}\E|\bar{Z}_s|^\alpha
  ds=2^\alpha e^{-\frac{\alpha t}{4}}\E|\bar{Z}_0|^\alpha$$
with $\E|\bar{Z}_0|^\alpha<+\infty$. The
third term of the right-hand-side is equal to
$e^{-\frac{\alpha t}{4}}\E\bigg|\int_{0}^{+\infty}e^{-\frac{s}{2}}K(\bar{Z}_s)ds\bigg|^\alpha$
where the expectation is finite according to Lemma
\ref{contmominv} below.

\end{adem}
\begin{alem}\label{contmominv}
 Let $(\zeta_t)_{t\geq 0}$ solve the SDE
$$d\zeta_t=2\sqrt{\nu}dB_t+cK(\zeta_t)dt-\frac{\zeta_t}{2}dt$$
for some real constant $c$.  Then,
$$\forall n\in\N^*,\;\exists C<+\infty,\;\forall \zeta_0,\;\E\left(\left(\int_0^{+\infty}
    \frac{e^{-\frac{s}{2}}}{|\zeta_s|}ds\right)^n\right)\leq C\mbox{
  and }\;\forall t\geq 0,\;\E\left(\left(\int_0^t
    \frac{ds}{|\zeta_s|}\right)^n\right)\leq C(1+t)^n.$$\end{alem}
\begin{adem}
The process ${R}_t\stackrel{\rm def}{=}|\zeta_t|^2$ solves the Cox-Ingersoll-Ross
stochastic differential equation
$$d{R}_t=4\sqrt{\nu {R}_t}d\beta_t+(8\nu-{R}_t)dt$$
with $d\beta_t=\frac{\zeta_t.dB_t}{\sqrt{{R}_t}}$ and does not
vanish since $8\nu=\frac{(4\sqrt{\nu})^2}{2}$. By the comparison
principle satisfied by this stochastic differential equation (see Theorem 3.7 p394 \cite{reyor}), for $t>s$,
$$\E\left(\frac{1}{\sqrt{R_t}}\bigg|R_s\right)\leq
\int_0^{+\infty}\frac{1}{\sqrt{r}}p(t-s,0,r)dr$$
where for $t>0$,
$p(t,0,r)=\frac{1}{8\nu(1-e^{-t})}e^{-\frac{r}{8\nu(1-e^{-t})}}$ denotes
the transition density from the state $0$. One deduces that
$\E\left(\frac{1}{\sqrt{R_t}}\bigg|R_s\right)\leq\sqrt{\frac{\pi}{8\nu(1-e^{s-t})}}$. By successive
conditionings, one deduces that
$$\forall n\in\N^*,\;\exists C_n>0,\;\forall
0=t_0<t_1<\hdots<t_n,\;\E\left(\frac{1}{\prod_{k=1}^n|\zeta_{t_k}|}\right)\leq
\frac{C_n}{\prod_{k=1}^n\sqrt{1-e^{t_{k-1}-t_k}}}.$$
As a consequence 
\begin{align*}
 \E\left(\left(\int_0^t
    \frac{ds}{|\zeta_s|}\right)^n\right)&=n!\int_0^t\int_{t_1}^t\hdots\int_{t_{n-1}}^t\E\left(\frac{1}{\prod_{k=1}^n|\zeta_{t_k}|}\right)dt_ndt_{n-1}\hdots dt_1\\
&\leq C\int_0^t\int_{t_1}^t\hdots\int_{t_{n-1}}^t\frac{dt_ndt_{n-1}\hdots dt_1}{\prod_{k=1}^n\sqrt{1-e^{t_{k-1}-t_k}}},
\end{align*}
Since $\forall s\in[0,t],\;\int_s^t\frac{1}{\sqrt{1-e^{s-r}}}dr\leq
\int_0^t\frac{1}{\sqrt{1-e^{-r}}}dr\leq C(1+t)$, one concludes that 
$\E\left(\left(\int_0^t
    \frac{ds}{|\zeta_s|}\right)^n\right)\leq C(1+t)^n$.
The finiteness of $\E\left(\left(\int_0^{+\infty}
    \frac{e^{-\frac{s}{2}}}{|\zeta_s|}ds\right)^n\right)$ is
obtained by a similar argument remarking that
$\forall s\geq
0,\;\int_s^{+\infty}\frac{e^{-\frac{r}{2}}}{\sqrt{1-e^{s-r}}}dr\leq \int_0^{+\infty}\frac{e^{-\frac{r}{2}}}{\sqrt{1-e^{-r}}}dr<+\infty$.

\end{adem}
\begin{alem}\label{contmomz}
  If for some
$\gamma\geq 2$,
$\E|X^1_0-X^2_0|^\gamma<+\infty$, then $\sup_{t\geq
  0}\E|Z_t|^\gamma<+\infty$.  
\end{alem}
\begin{adem}
By It\^o's formula, for $\delta\geq 2$, 
$$d|Z_t|^\delta=2\sqrt{\nu}\delta|Z_t|^{\delta-2}Z_t.dB_t+\left(2\nu\delta^2|Z_t|^{\delta-2}-\frac{\delta}{2}|Z_t|^\delta\right)dt.$$
For the choice $\delta=2$, since $\E|Z_0|^2<+\infty$, by a standard
localization argument we obtain that $d\E|Z_t|^2=(8\nu-\E|Z_t|^2)dt$
and deduce that $\sup_{t\geq
  0}\E|Z_t|^2<+\infty$. By induction, we then check that for all
$k\in\{1,\hdots,\lfloor\frac{\gamma}{2}\rfloor\}$, $d\E|Z_t|^{2k}=(8\nu
k^2\E|Z_t|^{2(k-1)}-k\E|Z_t|^{2k})dt$ so that 
$$\E|Z_t|^{2k}\leq e^{-kt}\E|Z_0|^{2k}+8\nu
k^2\int_0^te^{k(s-t)}\E|Z_s|^{2(k-1)}ds\mbox{ and }\sup_{t\geq
  0}\E|Z_t|^{2k}<+\infty.$$With the choice $\delta=\gamma$, we then
conclude that $d\E|Z_t|^{\gamma}=(2\nu
\gamma^2\E|Z_t|^{\gamma-2}-\frac{\gamma}{2}\E|Z_t|^{\gamma})dt$ and that $\sup_{t\geq
  0}\E|Z_t|^{\gamma}<+\infty$.
   
 \end{adem} 
\begin{alem}For $s\in [0,t/2]$, let ${\cal
    F}^t_s=\sigma(Z^{t/2}_t,(B_r-B_t)_{r\in[t-s,t]})$. Then
  $(\hat{B}^t_s\stackrel{\rm
    def}{=}B_{t-s}-B_t+\frac{1}{2\sqrt{\nu}}\int_{t-s}^tZ^{t/2}_rdr)_{s\in[0,t/2]}$ is a ${\cal F}^t_s$-Brownian motion and $(\hat{Z}^t_s\stackrel{\rm
    def}{=}Z^{t/2}_{t-s})_{s\in[0,t/2]}$ solves
$$d\hat{Z}^t_s=2\sqrt{\nu}d\hat{B}^t_s-aK(\hat{Z}^t_s)ds-\frac{1}{2}\hat{Z}^t_sds.$$
Moreover, $(Z^{t/2}_t,\int_{t/2}^te^{\frac{s-t}{2}}dB_s,\int_{t/2}^te^{\frac{s-t}{2}}K(Z^{t/2}_s)ds)$and $(\bar{Z}_0,\int_0^{t/2}e^{-\frac{s}{2}}(\bar{Z}_s/\sqrt{4\nu}\,ds-dB_s),\int_0^{t/2}e^{-\frac{s}{2}}K(\bar{Z}_s)ds)$ have the same distribution.
\label{timerev}
\end{alem}
\begin{adem}
One has $$\left(Z^{t/2}_t,\int_{t/2}^te^{\frac{s-t}{2}}dB_s,\int_{t/2}^te^{\frac{s-t}{2}}K(Z^{t/2}_s)ds\right)=\left(\hat{Z}^t_{0},\int_0^{t/2}e^{-\frac{s}{2}}(\hat{Z}^t_s/\sqrt{4\nu}\,ds-d\hat{B}^t_s),\int_0^{t/2}e^{-\frac{s}{2}}K(\hat{Z}^t_{s})ds\right).$$
Moreover, the trajectorial uniqueness for the stochastic differential equation \eqref{defzbar} implies that $(\bar{Z}_s,B_s)_{s\in[0,t/2]}$ and $(\hat{Z}^t_s,\hat{B}^t_s)_{s\in[0,t/2]}$ have the same distribution. Therefore the second statement is a consequence of the first one.

 Without the singularity of the Biot and Savart kernel at the origin,
 the first statement would be a consequence of \cite{Pa},  Theorem 2.2 and
 Corollary 2.4 (see also \cite{MNS} for more general results concerning
 the time-reversal of diffusion processes). To deal with this
 singularity, we use the smooth approximations $K_\varepsilon$ of this kernel defined at the end of the introduction. For $\varepsilon>0$, let $(Z^{t/2}_{\varepsilon,s})_{s\geq t/2}$ be the unique solution of the SDE
$$Z^{t/2}_{\varepsilon,s}=Z^{t/2}_{t/2}+2\sqrt{\nu}(B_s-B_{t/2})+\int_{t/2}^s
\left(aK_\varepsilon(Z^{t/2}_{\varepsilon,r})-\frac{1}{2}Z^{t/2}_{\varepsilon,r}\right)dr.$$
By an adaptation of the proof of the first statement in Proposition \ref{convdensnorm}, for all $s\geq t/2$, $Z^{t/2}_{\varepsilon,s}$ is
distributed according to the density $q_\infty$. By
Theorem 2.2 \cite{Pa} which deals with time-reversal of Brownian motions, one deduces that for $g:\R^2\rightarrow\R$ continuous and bounded and
$0\leq r\leq s\leq t/2$,
$$\E\left(\left(B_{t-s}-B_{t-r}+\frac{1}{2\sqrt{\nu}}\int_{t-s}^{t-r}
    Z^{t/2}_{\varepsilon,u}du\right)g(Z^{t/2}_{\varepsilon,t-r})\right)=0.$$
Hence
\begin{align}
   \bigg|\E\bigg((\hat{B}^t_s-\hat{B}^t_r)&g(Z^{t/2}_{t-r})\bigg)\bigg|\leq \frac{1}{2\sqrt{\nu}}\left|\E\left(g(Z^{t/2}_{\varepsilon,t-r})\int_{t-s}^{t-r}
         (Z^{t/2}_{u}-Z^{t/2}_{\varepsilon,u})du\right)\right|\notag\\&+\left|\E\left(\left(B_{t-s}-B_{t-r}+\frac{1}{2\sqrt{\nu}}\int_{t-s}^{t-r}
         Z^{t/2}_{u}du\right)(g(Z^{t/2}_{t-r})-g(Z^{t/2}_{\varepsilon,t-r}))\right)\right|
.\label{propmart}
\end{align}
Since $Z^{t/2}_u$ does not vanish, $\tau_\varepsilon\stackrel{\rm
  def}{=}\inf\{u\geq \frac{t}{2}:|Z^{t/2}_u|\leq \varepsilon\}$ goes to
infinity when $\varepsilon$ goes to $0$. The processes $(Z^{t/2}_u)_{u\in[t/2,t]}$ and
$(Z^{t/2}_{\varepsilon,u})_{u\in[t/2,t]}$ coincide on $\tau_\varepsilon\geq t$. By
Lebesgue's theorem, the second term of the right-hand-side of
\eqref{propmart} converges to $0$ as $\varepsilon\to 0$. So does the first
term since
\begin{align*}
   \left|\E\left(g(Z^{t/2}_{\varepsilon,t-r})\int_{t-s}^{t-r}
         (Z^{t/2}_{u}-Z^{t/2}_{\varepsilon,u})du\right)\right|&\leq
     \|g\|_\infty\int_{t-s}^{t-r}\E\left(1_{\{\tau_\varepsilon<t\}}\left|Z^{t/2}_{u}-Z^{t/2}_{\varepsilon,u}\right|\right)du\\
&\leq
     \|g\|_\infty\sqrt{\P(\tau_\varepsilon<t)}\int_{t-s}^{t-r}\sqrt{2\E\left(|Z^{t/2}_{u}|^2+|Z^{t/2}_{\varepsilon,u}|^2\right)}du\\&=2\|g\|_\infty(s-r)\sqrt{\P(\tau_\varepsilon<t)\E(|\bar{Z}_0|^2)}.
\end{align*}
Therefore $\E(\hat{B}^t_s-\hat{B}^t_r|Z^{t/2}_{t-r})=0$. 
Since ${\cal
  F}^t_r=\sigma(Z^{t/2}_{t-r},(B_u-B_t)_{u\in[t-r,t]})$ and
$\hat{B}^t_s-\hat{B}^t_r=B_{t-s}-B_{t-r}+\frac{1}{2\sqrt{\nu}}\int_{t-s}^{t-r}Z^{t/2}_udu$,
one deduces that $\E(\hat{B}^t_s-\hat{B}^t_r|{\cal F}^t_r)=0$.
Taking into
account the quadratic variation, one concludes that $\hat{B}_s$ is a
${\cal F}^t_s$-Brownian motion. The second statement follows easily.

\end{adem}

{\bf Acknowledgements.} We thank   two referees for  their careful reports on the  first version of this paper (in particular for  pointing out a typo in one of the  statements)  and for their remarks and questions that allowed us to   improve the presentation of  the results.

\end{document}